\documentclass[a4paper,11pt]{article}
\usepackage[utf8]{inputenc}
\usepackage{fullpage}
\usepackage{xargs}                      
\usepackage[pdftex,dvipsnames]{xcolor}  

\usepackage{fullpage,amsfonts,amsmath,amsthm,amssymb,mathtools}
\usepackage{tikz,color,algorithmic,algorithm,multirow,graphicx,caption,subcaption}
\theoremstyle{plain}
\newtheorem{theorem}{Theorem}
\newtheorem{corollary}[theorem]{Corollary}
\newtheorem{lemma}[theorem]{Lemma}
\newtheorem{remark}[theorem]{Remark}
\theoremstyle{definition}
\newtheorem{assumption}[theorem]{Assumption}
\newtheorem{definition}[theorem]{Definition}

\newcommand{\R}{\mathbf{R}}

\newcommand{\xs}{x^*}
\newcommand{\kk}{_k}%
\newcommand{\kz}{_0}%

\newcommand{\kpo}{_{k+1}}
\newcommand{\kmo}{_{k-1}}

\newcommand{\eqdef}{:=}

\newcommand{\mat}[1]{\begin{bmatrix}#1\end{bmatrix}}
\newcommand{\aln}{{\tilde{A}^{2\ell-1}}}
\newcommand{\apn}{{\tilde{A}^{2\ell + 1}}}
\newcommand{\al}{{\tilde{A}^{2\ell}}}
\newcommand{\g}{\tilde{g}}
\newcommand{\x}{x}
\newcommand{\Pre}{P^{-1}}
\newcommand{\PreT}{P^{-T}}

\usetikzlibrary{arrows}
\title{A Flexible Algorithmic Framework for Strictly Convex Quadratic Minimization}
\author{Liam MacDonald, Rua Murray and Rachael Tappenden}
\date{}

\begin{document}

\maketitle

\begin{abstract}
   This paper presents an algorithmic framework for the minimization of strictly convex quadratic functions. The framework is flexible and generic. At every iteration the search direction is a linear combination of the negative gradient, as well as (possibly) several other `sub-search' directions, where the user determines which, and how many, sub-search directions to include. Then, a step size along each sub-direction is generated in such a way that the gradient is minimized (with respect to a matrix norm), over the hyperplane specified by the user chosen search directions. Theoretical machinery is developed, which shows that any algorithm that fits into the generic framework is guaranteed to converge at a linear rate. Moreover, these theoretical results hold even when relaxation and/or symmetric preconditioning is employed. Several state-of-the-art algorithms fit into this scheme, including steepest descent and conjugate gradients.
\end{abstract}
\section{Introduction}

This work investigates an iterative scheme for solving the strongly convex quadratic optimization problem 
\begin{equation}\label{quadprob}
    \min\limits_{z \in \R^n} f(z) = \tfrac{1}{2}z^TAz - z^Tb,
\end{equation}
where $b \in \R^n$ and $A \in \R^{n \times n}$ is symmetric positive definite. Problem \eqref{quadprob} is equivalent to solving a system of linear equations, and by strong convexity, the unique solution is 
\begin{equation}\label{uniqueSol}
    z^* = A^{-1}b.
\end{equation}

Problem \eqref{quadprob} is the focus of much research, and this is justified due to its wide applicability in all areas of applied mathematics, 
economics, finance, and in natural and engineering sciences. Traditionally, direct methods that involve matrix factorizations have been used successfully to solve \eqref{quadprob}, although when $n$ is very large, such approaches can become computationally impractical. In a high-dimensional setting, iterative methods may be preferable, and such methods are the focus of this work.

Iterative methods for the solution of \eqref{quadprob} often make use of the gradient 
\begin{equation}\label{gradient}
    g(z) \eqdef \nabla f(z) = Az - b,
\end{equation}
(with the shorthand $g\kk = g(z\kk)$) as a search direction. Arguably, the best known example of a gradient based method for solving \eqref{quadprob} is the Steepest Descent (SD) method \cite{cauchy1847methode}, developed by Cauchy in 1847, which uses the negative gradient as the search direction and uses a step size that exactly minimizes the function value along that direction. The convergence properties of SD have been studied in \cite{Akaike1959,Nocedal2002,forsythe1951Asymptoticproperties}, and those works show that the method converges at a linear rate, and that asymptotically, the gradient oscillates in a subspace spanned by the left- and right-most eigenvectors. SD has been extended to a family methods, all of which use the gradient as the search direction, but choose a step size that exactly minimizes the norm of the gradient, where the norm is defined in terms of the Hessian, $A$, \cite{krasnosel1952iteration, Pronzato_2005, macdonald24}. 

While the aforementioned methods are all guaranteed to solve \eqref{quadprob}, the oscillating/`zig-zagging' nature of the successive gradients can lead to poor performance in practice, particularly on ill-conditioned problems. Previously, it was thought that using the gradient as the search direction that was the cause of the slow computational behaviour, but the works \cite{Raydan2002,TwoPointBarzilaiBorwein} show that keeping the gradient as the search direction but adjusting the step size (in \cite{Raydan2002} using relaxation/damping, and in \cite{TwoPointBarzilaiBorwein} using a novel step-size based on approximating the Hessian by a multiple of the identity) can lead to noticable improvements in practice. Other works that have considered modifications to the step size with the aim of improving the practical convergence properties of gradient based methods include \cite{dai2003alternate,Huang2021BBTermination,zhou2006gradient,DaiYuanAnalysis}.

Another approach to enhancing the practical performance of algorithms for \eqref{quadprob} is to adapt the search direction, possibly by using a combination of several `sub'-directions in the iterative update. An example of an algorithm that falls into this category is the momentum/heavy ball algorithm of Polyak \cite{polyak1964some}, where the current point is updated by taking the sum of a step along the negative gradient and a step along the direction formed from the difference between the past two iterates. Polyak's momentum was extended by Nesterov \cite{Nesterov1983}, who developed the accelerated gradient method. Acceleration can be viewed as applying a correction to the momentum term, to further support practical performance. These methods are the basis of many state-of-the-art algorithms in machine learning, and remain an active area of research \cite{pmlr-v202-wang23l, goujaud2022quadratic, van2017fastest}.

Conjugate Gradients (CG) \cite{Hestenes1952,nazareth1977conjugate} and the Conjugate Residual (CR) method both generate search directions that are conjugate, and where each direction can be, equivalently, expressed as a linear combination of the negative gradient and the previous search direction (although CG and CR each use a different norm to compute the coefficients in the iterative update). The extension of the conjugate gradient method to different norms is difficult to find in the literature, as has been recently noted \cite{lim2024conjugate, goujaud2022quadratic} and independently reproduced, but results exist in several earlier papers \cite{ashby1990taxonomy, broyden1996new, fasano2015framework}. Moreover, the series of papers \cite{brezinski1997projection, brezinski1996variations,brezinski1999multiparameter} considers modifying the search direction via a projection based framework, and by enforcing a conjugacy type condition on the search direction, CG is recovered.

Another algorithm relevant to the current work is Forsythe's $s$-directional methods \cite{Forsythe1968}, where at each iteration the search direction is formed as a linear combination of the sub-directions $g\kk, Ag\kk,\dots, A^{s-1}g\kk$, for some fixed $s\geq 1$. These methods are typically not used in practice because they are expensive (requiring multiple matrix vector products at every iteration), and they often display oscillating behaviour, similar to that of steepest descent. Nevertheless, \cite{Pronzato2009} shows that using $s=1$ and $s=2$ on alternate iterations can break the zig-zagging behaviour. The original algorithm is also shown to be equivalent to taking $s$ steps of the conjugate gradient method and then restarting \cite{luenberger2008linear, Pronzato2009}.

The Barzilai-Borwein (BB) method \cite{TwoPointBarzilaiBorwein} is another algorithm that can be applied to solve~\eqref{quadprob}. At each iteration, the BB method uses the negative gradient as the search direction, and the formula for the step size involves differences of successive past gradients, and successive past iterates. Its popularity persists because the iterates are inexpensive and it has similar practical performance to CG on certain problem instances. An extension of the BB method, also for solving \eqref{quadprob}, is Limited Memory Steepest Descent (LMSD) \cite{fletcher2012limited, curtis2018linear}, where, at every iteration, LMSD also uses the negative gradient as the search direction. However, for LMSD, the step sizes are calculated (and then used) in `batches of $m$', and are based upon the Ritz values of an $m\times m$ matrix that is constructed using the previous $m$ gradients $g_{k-m},\dots,g_{k-1}$. (When $m=1$, the BB method is recovered). While the BB and LMSD methods can be used to solve \eqref{quadprob}, they fall outside the scope of the current work, because they use a single search direction at every iteration (the negative gradient), and they use step sizes that are generated using `approximations', while this work restricts its attention to exact step sizes.

The goal of this work is to suggest a generic scheme for solving problem \eqref{quadprob}, where at each iteration the search direction is made up of a linear combination of sub-search directions, where the user determines which, and how many, sub-search directions to use. The scheme is flexible. As long as one of the sub-search directions is the negative gradient, the number of sub-search directions, and the directions themselves, is arbitrary, and can change from one iteration to the next. A formula is provided for the step-sizes, and it is formulated such that the norm of the gradient is minimized over the hyperplane defined by the user-chosen sub-search directions. Theoretical convergence results are provided, which show that any algorithm fitting the generic scheme is guaranteed to converge at a linear rate. Hence, this work also provides a theoretical machinery for future algorithms: any new algorithm that fits the generic scheme described in this work is automatically equipped with the stated convergence guarantees.

The aforementioned works are known to converge linearly, with the rate typically depending on the condition number. Problem conditioning also impacts the practical performance of numerical methods, which tend to behave poorly on ill-conditioned problem instances. Preconditioning is a technique that can be employed, with the goal of adjusting the spectrum of the iteration matrix, thereby improving problem conditioning and algorithm performance. Jacobi, Gauss-Seidel, and SSOR are all algorithms that can be expressed as preconditioned gradient descent methods for problem \eqref{quadprob} and the development of effective preconditioners is an active area of research, see for example \cite{jahani2021iclr,sadievTappenden2022stochastic,maddison2021}. The current work allows for the inclusion of (symmetric) preconditioning, with all theoretical results also holding if preconditioning is employed.

\subsection{Notation and Preliminaries}

The following assumption and preliminaries are utilised in this work.

\begin{assumption}\label{Assume1}
The matrix $A\in \R^{n\times n}$ in \eqref{quadprob} is symmetric and positive definite.
\end{assumption}

By Assumption \ref{Assume1}, \eqref{quadprob} is strongly convex. Moreover, because $A$ is positive definite, it has a unique, positive definite square root, denoted by $A^{1/2}$. The convention $A^0 = I$ is adopted, the eigenvalues of $A$ are denoted by $0<\lambda_n\leq \cdots \leq \lambda_1$, the condition number is $\kappa = \lambda_1/\lambda_n$, and $e_j\in \R^n$ denotes the $j$th column of the $n\times n$ identity matrix.

\begin{definition}
Let $B \in \R^{n \times n}$ be a symmetric positive definite matrix. Then, for any $v \in \R^n$, define the norm $\|v\|_{B}=\sqrt{v^TBv}$, and the inner product $ \langle v, v \rangle_B = \langle v, Bv \rangle = \langle Bv, v \rangle$.
\end{definition}

\subsection{Paper outline}

The paper is organized as follows. Section~\ref{sec:genericscheme} introduces the algorithmic framework for minimizing strictly convex quadratics. Theoretical results are presented showing that any algorithm fitting the framework is guaranteed to converge at a linear rate, even when relaxation and symmetric preconditioning are used. Section~\ref{sec:Algsfittingframework} shows that several existing state-of-the-art methods fit the generic scheme for \eqref{quadprob}, including steepest descent, conjugate gradients, and Forsythe's $s$-gradient methods. Numerical experiments are presented in Section~\ref{sec:NumericalExp}, and a discussion and conclusion is given in Section~\ref{sec:Conclusion}.

\section{A flexible scheme for quadratic minimization}
\label{sec:genericscheme}

This section describes a flexible scheme for solving \eqref{quadprob}. The scheme is `multi'-directional, allowing for several `sub'-directions to make up the search direction, and the step size is chosen to exactly minimize the norm (where the norm is possibly defined in terms of the (preconditioned) Hessian), of the gradient over the hyperplane defined by the search directions. In particular, given an initial point $z\kz \in \R^n$, for all $k\geq 0$ let $1\leq m\kk \leq n$ and consider iterates of the form 
\begin{equation}\label{xupdate}
    z\kpo = z\kk - \alpha\kk^{(1)}w\kk^{(1)} - \cdots -\alpha\kk^{(m\kk)}w\kk^{(m\kk)} = z\kk - W\kk a\kk,
\end{equation}

where  
\begin{equation}\label{W}
    W\kk \eqdef \mat{w\kk^{(1)}, & \cdots, & w\kk^{(m\kk)}} \in \R^{n\times m\kk},
\end{equation} 
\begin{equation}\label{a}
    a\kk \eqdef \mat{\alpha\kk^{(1)},& \cdots, & \alpha\kk^{(m\kk)}}^T\in \R^{m\kk}.
\end{equation}
The columns of $W\kk$ in \eqref{W} are the $m\kk$ `sub'-search directions, and the elements of $a\kk$ are the corresponding `sub'-step sizes \eqref{a}. The parameter $m\kk$ is iteration dependent, so the number of sub-search directions is allowed to vary from one iteration to the next. 

As has been mentioned, iterative schemes for solving \eqref{quadprob} can benefit from preconditioning, and the generic scheme described in this work allows for the possibility of a preconditioner. Recall that the solution to \eqref{quadprob} is equivalent to the solution of a system of linear equations, so letting $\Pre \in \R^{n\times n}$ be full rank, one can consider the equivalent problem of solving the symmetric preconditioned system:
\begin{eqnarray}\label{precondf}
 \Pre A \PreT x = \Pre b,\quad \text{where}\quad \PreT x = z.
\end{eqnarray}

In this work preconditioned variables are typically denoted using a `tilde'. For convenience, define
\begin{equation}\label{PA}
    \tilde{A} = \Pre A \PreT,\qquad \text{and} \qquad \tilde{b} = \Pre b,
\end{equation} 
so that the preconditioned gradient is denoted by
\begin{equation}\label{Pgradient}
    \g\kk = \Pre (A \PreT \x\kk - b) = \Pre g(z\kk),
\end{equation}

Preconditioning can be incorporated into \eqref{xupdate}, giving the update 
\begin{equation}\label{PreconNoRelax}
    \x\kpo = \x\kk - \Pre W\kk a\kk,
\end{equation}

The following assumptions are made regarding the vectors $w\kk^{(1)},\dots,w\kk^{(m\kk)}$ and the preconditioner $\Pre$.
\begin{assumption}\label{Assumepi}
For all $k\geq 0$, suppose that $W\kk$ in \eqref{W} satisfies the following properties:
\begin{itemize}
    \item[(i)] $W\kk$ has full rank; and
    \item[(ii)] $g(z\kk) \in {\rm col}(W\kk)$.
\end{itemize}   
i.e., $w\kk^{(1)},\dots,w\kk^{(m\kk)} \in \R^{n}$ are linearly independent, and $g(z\kk) \in {\rm span}\{w\kk^{(1)},\dots,w\kk^{(m\kk)}\}$.
\end{assumption}
\begin{assumption}\label{Assume_Pre}
    The preconditioner $\Pre\in \R^{n\times n}$ is nonsingular.
\end{assumption}
\begin{remark}
    If Assumptions \ref{Assume1} and \ref{Assume_Pre} hold, then $\tilde A$ in \eqref{PA} is symmetric and positive definite.
\end{remark}

Given search directions $w\kk^{(i)}$ for $1\leq i\leq m\kk$ and a preconditioner $\Pre$, at any iteration $k\geq 0$, the step size used in the generic scheme is given in the following lemma.

\begin{lemma}\label{stepsizeminimizegrad}
Let $f$ be defined in \eqref{quadprob}, let Assumptions~\ref{Assume1}, \ref{Assumepi} and \ref{Assume_Pre} hold, and fix $\ell \in \{0,\tfrac12,1,\tfrac32,2,\tfrac52,3,\dots\}$. Given $\x\kz\in \R^n$, let the iterates be as in \eqref{PreconNoRelax}, where
\begin{equation}\label{GeneralisedStepsize}
    a\kk \eqdef (W\kk^T \PreT \apn \Pre  W\kk)^{-1}W\kk^T \PreT \al \g\kk.
\end{equation}
Then 
\begin{equation}\label{mina2lnorm}
    a\kk = \arg\min\limits_{a} \|\g\kpo\|_{\aln}^2.
\end{equation}
\end{lemma}
\begin{proof}
By Assumptions~\ref{Assume1}, \ref{Assumepi} and \ref{Assume_Pre}, $W\kk$ has full rank and $W\kk^T \apn W\kk$ is nonsingular. Then
\begin{align}
    \|\g\kpo\|_{\aln}^2 
    &\overset{\eqref{Pgradient}+\eqref{PreconNoRelax}}{=} \| \Pre A \PreT (x\kk - \Pre W\kk a)-\Pre b\|_{\aln}^2 \notag\\
    &\overset{\eqref{Pgradient}+\eqref{PA}}{=} \| \g\kk - \tilde{A} \Pre W\kk a\|_{\aln}^2 \label{gkpovgk}\\
    &\overset{\phantom{\eqref{Pgradient}+\eqref{PA}}}{=} \|\g\kk \|_{\aln}^2 - 2 \langle \Pre W\kk a,\g\kk \rangle_{\al} + \|  \Pre W\kk a \|_{\apn}^2\notag.
\end{align}
Differentiating w.r.t. $a$, setting to zero, and rearranging, gives \eqref{GeneralisedStepsize}.
\end{proof}

\begin{remark}
    Determining the step size \eqref{GeneralisedStepsize} involves solving an $m\kk\times m\kk$ system. Typically $m\kk$ is small (usually $m\kk\ll n$), so this is inexpensive.
\end{remark}

The Flexible Algorithmic Framework (which, henceforth, will also be referred to as the Flexible Scheme) for solving \eqref{quadprob} is presented now.

\begin{definition}[Flexible Algorithmic Framework for \eqref{quadprob}]\label{def:Algorithms}
   Let $f$ be given in \eqref{quadprob}, let Assumptions \ref{Assume1}, ~\ref{Assumepi}~and~\ref{Assume_Pre} hold, and fix $\ell \in \{0,\tfrac12,1,\tfrac32,2,\tfrac52,3,\dots\}$ and $\omega \in (0,2)$. Given an initial point $\x\kz \in \R^n$, for all $k\geq 0$, the iterations of a Flexible Scheme for \eqref{quadprob} are defined by 
   \begin{equation}\label{relaxedIterates}
    \x\kpo = \x\kk - \omega \Pre W\kk a\kk,
   \end{equation}
   where $W\kk$ and $a\kk$ are given in \eqref{W} and \eqref{a}, respectively.
\end{definition}

The update \eqref{relaxedIterates}  in Definition~\ref{def:Algorithms} is slightly different from \eqref{PreconNoRelax}; it involves a relaxation or damping parameter $\omega \in (0, 2)$. It is known that relaxation can improve the practical behaviour of certain iterative methods, which motivates its inclusion here. If $\omega = 1$ (no relaxation), then \eqref{relaxedIterates} is equivalent to \eqref{PreconNoRelax}, and the update \eqref{xupdate} is recovered when the identity matrix is chosen as the preconditioner, i.e., $\Pre=I$. 

\paragraph{Special case when $W\kk$ is nonsingular} If $W\kk$ is nonsingular (i.e., $m\kk = n$ and $W\kk$ has full rank), then $a\kk \overset{\eqref{GeneralisedStepsize}}{=} W\kk^{-1}P P^T A^{-1}g(z\kk).$
Hence,
\begin{equation*}
    x\kpo \overset{\eqref{PreconNoRelax}}{=} x\kk - \Pre W\kk a\kk 
    = x\kk -  P^T A^{-1} g(z\kk) 
    \overset{\eqref{gradient}}{=} P^T A^{-1}b,
\end{equation*}
so that the change of variables returns $\PreT x\kpo = z^*$,
i.e., for \emph{any} nonsingular $W\kk$, convergence is achieved in exactly one step, as expected.

\subsection{Theoretical properties}\label{sec:theory}

This section presents theoretical results for the Flexible Scheme given in Definition~\ref{def:Algorithms}. The convergence proof proceeds by showing that, for fixed $\ell$, the Flexible Scheme performs at least as well as a preconditioned variant of the Relaxed $\ell$-MGD method in \cite{macdonald24}. Because these methods are a key component of the theoretical analysis, they are summarised here for convenience. 

The Relaxed $\ell$-MGD methods in \cite{macdonald24} are a special case of Definition~\ref{def:Algorithms} when $\Pre = I$ and $m\kk = 1$ for all $k\geq 0$ (i.e., they satisfy Assumption~\ref{Assumepi}, so $w\kk^{(1)} = g(z\kk)= g(x\kk)$ for all $k\geq 0$); see also Section~\ref{sec:relaxedellMG}. Substituting these parameter choices into \eqref{GeneralisedStepsize} gives $a\kk = (g\kk^TA^{2\ell}g\kk)/(g\kk^TA^{2\ell+1}g\kk)\in \R$, which matches the step length given in \cite{macdonald24}. Because this work allows preconditioning, the Relaxed $\ell$-MGD methods are now modified to give a family of Preconditioned Relaxed $\ell$-MGD methods (that also fit Definition~\ref{def:Algorithms}). In particular, let $\Pre$ satisfy Assumption~\ref{Assume_Pre}, and let $m\kk = 1$ for all $k\geq 0$, with $w\kk^{(1)} = g(z\kk)$ for all $k\geq 0$ to ensure Assumption~\ref{Assumepi} holds. It is convenient to define
\begin{equation}\label{familystep}
    \theta\kk \eqdef \arg\min\limits_{\theta} \|\g\kk - \theta \tilde{A} \g\kk\|^2_{\aln}
    = \frac{\langle \g\kk, \g\kk \rangle_\al }{\|\g\kk\|_{\apn}^2},
\end{equation}
and notice that $\theta\kk \equiv a\kk \in \R$ when the previously mentioned algorithm parameters are substituted into \eqref{GeneralisedStepsize} (and hence the minimization property \eqref{mina2lnorm} also holds for this $\theta\kk$). 

Now, consider the following. By Assumption~\ref{Assumepi}, $g(z\kk) \in {\rm span}\{w\kk^{(1)},\dots,w\kk^{(m)}\}$, so without loss of generality, for any $k\geq 0$ let $g(z\kk)$ be the first column of $W\kk$, i.e., $w\kk^{(1)} = g(z\kk)$. Then, the choice $a\kk = \theta\kk e_1 = [\theta\kk, 0, 0, \dots, 0]^T$ recovers the Preconditioned Relaxed $\ell$-MGD method. The proof of convergence of the Flexible Scheme in Definition~\ref{def:Algorithms} proceeds by establishing that, due to the minimization property \eqref{mina2lnorm}, computing $a\kk$ as in \eqref{GeneralisedStepsize} leads to a reduction in the norm of the preconditioned gradient that is at least as large as that achieved by taking $a\kk = \theta\kk e_1$. This implies that the Flexible Algorithmic Framework inherits all of the convergence properties of the Preconditioned Relaxed $\ell$-MGD method, and this argument is formalised now.

\begin{lemma}\label{lem:atleastMGreduction}
Let $f$ be given in \eqref{quadprob}, let Assumptions~\ref{Assume1},~\ref{Assumepi}~and~\ref{Assume_Pre} hold, fix $\ell \in \{0,\tfrac12,1,\tfrac32,2,\tfrac52,3,\dots\}$ and $\omega \in (0,2)$, and let $\theta\kk$ be as in~\eqref{familystep}. Given an initial point $\x\kz \in \R^n$, for all $k\geq 0$, let the iterates of the Flexible Scheme be given in Definition~\ref{def:Algorithms} where $W\kk$ and $a\kk$ are given in \eqref{W} and \eqref{a}, respectively. Then, 
    \begin{equation}\label{eqn:atleastMGD}
\|\g\kk -\omega \tilde{A} \Pre W\kk a\kk \|^2_\aln \leq \|\g\kk - \omega \theta\kk \tilde{A} \g\kk\|^2_\aln. 
\end{equation}
\end{lemma}
\begin{proof}
If $m\kk = 1$ then by Assumption~\ref{Assumepi}, $W\kk = g(z\kk)$, so $\x\kk -\omega \Pre W\kk a\kk = \x\kk - \omega \theta\kk \Pre g(z\kk)$. Multiplying through by $\tilde A$ and subtracting $\tilde b$ from both sides, shows that \eqref{eqn:atleastMGD} holds with equality. Hence, for the remainder of this proof, suppose that $m\kk \geq 2$.

By Assumption~\ref{Assumepi}(i), $w\kk^{(1)},\dots,w\kk^{(m\kk)}$ are linearly independent, and by Assumption~\ref{Assume_Pre}, so too are $\Pre w\kk^{(1)},\dots,\Pre w\kk^{(m\kk)}$. Hence, w.l.o.g. assume that they form an orthonormal set in $\R^n$ w.r.t. the inner product induced by $ \apn $, i.e., let $\Pre w\kk^{(1)},\dots, \Pre w\kk^{(m\kk)}$ be $\apn $--conjugate (i.e., applying Gram-Schmidt achieves this).
Thus, 
\begin{equation}\label{OrthogonalP}
    W\kk^T \PreT \apn \Pre W\kk = I,
\end{equation}
which can be substituted into \eqref{GeneralisedStepsize}, giving
\begin{equation} \label{OrthogonalSteps}
    a\kk = W\kk^T \PreT \al \g\kk.
\end{equation}
By Assumption~\ref{Assumepi}(ii), $g(z\kk) \in {\rm col}(W\kk)$, so w.l.o.g setting $\Pre w\kk^{(1)} =  \g\kk/\|\g\kk\|_\apn$ (where $\Pre w\kk^{(1)}$ is a unit vector because of the previous orthonormalization process) gives
\begin{equation}\label{alphaLforms}
    \theta\kk \overset{\eqref{familystep}}{=} \frac{\langle \g\kk, \g\kk \rangle_\al }{\|\g\kk\|_{\apn}^2} = \frac{\langle \Pre w\kk^{(1)}, \g\kk\rangle_{\al}}{\|\g\kk\|_\apn},
\end{equation}
and
\begin{equation}\label{OrthogonalAlphaNice}
    \alpha\kk^{(1)} \overset{\eqref{OrthogonalSteps}}{=} \langle \Pre w\kk^{(1)}, \g\kk\rangle_{\al} \overset{\eqref{alphaLforms}}{=} \theta\kk \|\g\kk\|_\apn.
\end{equation}
Note that for $\omega \in (0, 2)$, $0 \leq \omega(2 - \omega) \leq 1 $. Then,
\begin{eqnarray}
\| \g\kk - \omega \tilde{A} \Pre W\kk a\kk \|^2_\aln  
    &=& \|\g\kk\|_{\aln}^2  - 2\omega \langle \Pre W\kk a\kk,\g\kk \rangle_{\al}  + \omega^2\|  \Pre W\kk a\kk\|_{\apn}^2  \notag \\
    &\overset{\eqref{OrthogonalP}+\eqref{OrthogonalSteps}}{=}& \|\g\kk\|_{\aln}^2 - 2\omega \| a\kk\|_2^2 + \omega^2\| a\kk\|_2^2 \notag\\
    &\overset{\eqref{a}+\eqref{OrthogonalAlphaNice}}{=}& \|\g\kk\|_{\aln}^2 -2\omega \theta\kk^2 \|\g\kk\|^2_\apn+ \omega^2 \theta\kk ^2 \|\g\kk\|^2_\apn\notag\\
    && - \omega(2-\omega) \sum_{i= 2}^{m\kk} (\alpha\kk^{(i)})^2\notag\\
    &\overset{\eqref{alphaLforms}}{\leq}& \|\g\kk\|_{\aln}^2 -2\omega \theta\kk \langle \g\kk, \g\kk \rangle_{\al}+ \omega^2 \theta\kk ^2 \|\g\kk\|^2_\apn\notag\\
    &=& \|\g\kk - \omega \theta\kk \tilde{A} \g\kk\|^2_\aln. \notag \label{RelaxedFamilyNormDef}
\end{eqnarray}
\end{proof}

In words, Lemma~\ref{lem:atleastMGreduction} shows that, because of the minimization property \eqref{mina2lnorm}, the gradient norm due to the step length generated by \eqref{GeneralisedStepsize} will be no more than that generated by the step length \eqref{familystep} used for the Preconditioned Relaxed $\ell$-MGD method \cite{macdonald24}.  Now, Lemma~\ref{lem:atleastMGreduction} can be used to establish convergence of the Flexible Scheme. Let the eigenvalues of $\tilde A$ be denoted by $0 < \tilde{\lambda}_n \leq \cdots \leq \tilde{\lambda}_1$, and define $\tilde \kappa = \tilde \lambda_1/\tilde \lambda_n$ to be the condition number of $\tilde A$.
Define
\begin{eqnarray}\label{gmonodecreaserelaxedc}
c(\omega) \eqdef \left(1 - \omega(2-\omega)\frac{4\tilde{\kappa}}{(\tilde{\kappa} + 1)^2}\right),
\end{eqnarray}
and note that $0<c(\omega)<1$ for $\omega \in (0,1)$.  Theorem~5 in \cite{macdonald24} shows that (with minor modifications to allow for preconditioning):
\begin{equation}\label{relaxedgraddecrease}
\|\g\kk - \omega \theta\kk \tilde{A} \g\kk\|_\aln^2  \leq c(\omega)\|\g\kk\|_\aln^2 .
\end{equation}

Now, \eqref{relaxedgraddecrease} can be combined with Lemma~\ref{lem:atleastMGreduction} to give the following result, which establishes convergence of the Flexible Scheme in Definition~\ref{def:Algorithms}.
\begin{theorem}
\label{thm:linearconvergence}
    Let $f$ be given in \eqref{quadprob}, let Assumptions~\ref{Assume1},~\ref{Assumepi}~and~\ref{Assume_Pre} hold, fix $\ell \in \{0,\tfrac12,1,\tfrac32,2,\tfrac52,3,\dots\}$ and $\omega \in (0,2)$, and let $\theta\kk$ be as in~\eqref{familystep}. Given an initial point $\x\kz \in \R^n$, for all $k\geq 0$, let the iterates of the Flexible Scheme be given in Definition~\ref{def:Algorithms} where $W\kk$ and $a\kk$ are given in \eqref{W} and \eqref{a}, respectively. Then,
\begin{align}\label{eqn:gradientlinear}
    \|\g\kpo\|^2_{\aln} 
    \leq (c(\omega))^k\|\g\kz\|_\aln^2,
\end{align}
\end{theorem}
\begin{proof}
    \begin{eqnarray*}
\|\g\kpo\|^2_\aln
    &=& \|\g\kk - \omega \tilde{A} \Pre W\kk a\kk\|^2_\aln \notag \\
    &\overset{\eqref{eqn:atleastMGD}}{\leq}& \|\g\kk - \omega \theta\kk \tilde{A} \g\kk\|^2_\aln \\
    &\overset{\eqref{relaxedgraddecrease}}{\leq}& c(\omega)\|\g\kk\|_\aln^2.\notag
\end{eqnarray*}
Finally, applying a recursion argument gives the result. 
\end{proof}

\begin{corollary}
    Let the conditions of Theorem~\ref{thm:linearconvergence} hold. If $\omega=1$, then
    \begin{equation}
        \|\g\kpo\|_\aln^2 \leq \left(\frac{\tilde{\kappa} - 1}{\tilde{\kappa} + 1}\right)^2\|\g\kk\|^2_\aln .
    \end{equation}
\end{corollary}

It is possible to extend these results and show that for any algorithm that fits the Flexible Algorithmic Framework in Definition~\ref{def:Algorithms}, (1) the distance of the current iterate from optimality converges linearly with respect to an appropriately chosen norm; (2) the function values converge linearly; and (3) an explicit iteration complexity result is available for every algorithm.
\begin{lemma}
Let the conditions of Theorem~\ref{thm:linearconvergence} hold, and let $\epsilon >0$. Then
\begin{enumerate}
    \item $\|\x\kk-\xs\|_\apn^2 \leq (c(\omega))^k \; \|\x\kz-\xs\|_\apn^2$;
    \item $f(\x\kk)-f^* \leq \tilde{\kappa}^{2\ell}(c(\omega))^{k} \, (f(\x\kz)-f^*)$;
    \item and $f(\x_{K})-f^* \leq \epsilon$, where
\begin{equation}\label{eq:Khat}
    K > \frac{1}{\omega(2-\omega)}\frac{(\tilde{\kappa}+1)^2}{4\tilde{\kappa}}\left(\ln{\frac{\tilde{\kappa}^{2\ell}\left(f(\x\kz)-f^*\right)}{\epsilon}} \right).
\end{equation}
\end{enumerate}
\end{lemma}
\begin{proof}
Notice that $\|\x\kk-\xs\|_\apn^2 = \|\tilde{A}(\x\kk-\xs)\|_\aln^2 = \|\g\kk\|_\aln^2,$ so 1. holds. The proofs of 2. and 3. follow immediately from the proofs of Theorems 10(2) and 13(2) in \cite{macdonald24}.
\end{proof}

\section{Algorithms that fit the flexible scheme}\label{sec:Algsfittingframework}
The scheme described in Definition~\ref{def:Algorithms} is deliberately general so as to capture as wide a range of algorithm designs as possible. Here, it is shown that several well-known existing algorithms fit within this scheme.

\subsection{Relaxed $\ell$-Minimal Gradient Descent Methods}\label{sec:relaxedellMG}
Recall Section~\ref{sec:theory}, which shows that the Relaxed $\ell$-MGD methods \cite{macdonald24} fit Definition~\ref{def:Algorithms} used with $P=I$, and $m\kk = 1$ for all $k\geq 0$ so $W\kk = w\kk^{(1)} = g(z\kk)$.  The well known Steepest Descent \cite{Akaike1959}, and the Minimal Gradient \cite{krasnosel1952iteration} algorithms are contained within this family, so by extension, they are also members of the Flexible Scheme (Definition~\ref{def:Algorithms}). In particular, using the algorithm parameters just mentioned, and fixing $\ell =0$ gives the Steepest Descent method (i.e., at every iteration the step size $a\kk = \theta\kk$ is that which minimizes the $A^{-1}$-norm of the gradient), while the choice $\ell = 1/2$ gives the Minimal Gradient method (i.e., where the step size $a\kk = \theta\kk$ is that which minimizes the $2$-norm of the gradient). Convergence of these methods  is established in \cite{macdonald24} when $P=I$ (i.e., without preconditioning). 

Importantly, the theoretical results in this work extend those of  \cite{macdonald24} and show that this family of algorithms is guaranteed to converge when symmetric preconditioning with a nonsingular preconditioner $\Pre$ is employed, giving a family of Preconditioned Relaxed $\ell$-MGD methods.

\subsection{Conjugate Direction methods}\label{sec:CG}

The (Preconditioned) Conjugate Gradient (CG) algorithm is a state-of-the-art algorithm for solving problem \eqref{quadprob}. For simplicity, suppose that $P=I$ (i.e., no preconditioning) and note that the choice $\omega = 1$ is appropriate for CG. Given an initial point $x\kz$, and defining $d_{-1} = 0 \in \R^n$ and $\beta_{-1}^{CG}=0\in \R$, for all $k\geq 0$, the iterates of CG are generated as
\begin{equation}\label{xdCG}
x\kpo = x\kk + \alpha\kk^{CG} d\kk, \quad \text{where}\quad d\kk = -g\kk + \beta\kmo^{CG} d\kmo.
\end{equation}
The parameter $\beta\kmo^{CG}$ ensures that the directions generated by the algorithm \eqref{xdCG} are conjugate, while the parameter $\alpha\kk^{CG}$ is chosen in such a way that the function value is minimized along the search direction, i.e., $\alpha\kk^{CG} = \arg \min_{\alpha} f(x\kk + \alpha d\kk)$ ($\equiv \arg \min_{\alpha} \|\nabla f(x\kk + \alpha d\kk)\|_{A^{-1}}^2$). It is straightforward to show that defining
\begin{equation}\label{alphabeta}
    \alpha\kk^{CG} = -\frac{d\kk^Tg\kk}{d\kk^TAd\kk}\qquad \text{and} \qquad \beta\kmo^{CG} = \frac{d\kmo^TAg\kk}{d\kmo^TAd\kmo},
\end{equation}
give parameters that satisfy these conditions.

Now, consider setting $W\kk = \mat{g\kk,& s\kk}$ where $s\kk = x\kk - x\kmo$ in the Flexible Scheme in Definition~\ref{def:Algorithms}.
Substituting this $W\kk$ into \eqref{GeneralisedStepsize} (and fixing $\ell=0$) gives the step sizes
\begin{eqnarray}\label{FlexSchemeCGsteps}
    a\kk = \mat{g\kk^TAg\kk & g\kk^TAs\kk\\ s\kk^TAg\kk & s\kk^TAs\kk}^{-1}\mat{g\kk^Tg\kk\\s\kk^Tg\kk} =
   \frac1{\Delta}\mat{s\kk^TA s\kk g\kk^Tg\kk - g\kk^TAs\kk s\kk^Tg\kk \\
    -g\kk^TAs\kk g\kk^Tg\kk + g\kk^T A g\kk s\kk^T g\kk } = \mat{\alpha\kk^{(1)} \\ \alpha\kk^{(2)}},
\end{eqnarray}
where 
\begin{equation}
    \Delta = g\kk^TAg\kk s\kk^TAs\kk - \left(g\kk^TAs\kk\right)^2.
\end{equation}
The goal now is to show that the iterates generated by the CG algorithm are equivalent to those of the Flexible Scheme with the specified search directions in $W\kk$. To see this, note that the CG algorithm can be written in terms of $g\kk$ and $s\kk$ via
\begin{eqnarray}
    d\kmo &=& \tfrac{1}{\alpha\kmo^{CG}} \left( x\kk - x\kmo \right) = \tfrac{1}{\alpha\kmo^{CG}} s\kk,\label{AltDkmo}\\
    \beta\kmo^{CG} &=& \frac{d\kmo^T A g\kk}{d\kmo^T A d\kmo} \overset{(\ref{AltDkmo})}{=} \frac{\tfrac{1}{\alpha\kmo^{CG}} s\kk^TAg\kk}{\tfrac{1}{\left(\alpha\kmo^{CG}\right)^2}s\kk^TAs\kk} = \alpha\kmo^{CG} \frac{s\kk^TAg\kk}{s\kk^TAs\kk}\label{AltBetakmo}.
\end{eqnarray}
Substituting \eqref{AltDkmo} and \eqref{AltBetakmo} into \eqref{xdCG} gives
\begin{eqnarray}\label{xkCGgksk}
    x\kpo &\overset{(\ref{xdCG}+\ref{AltDkmo})}{=}& x\kk -\alpha\kk^{CG} g\kk + \alpha\kk^{CG} \frac{\beta\kmo^{CG}}{\alpha\kmo^{CG}}s\kk\notag\\
    &\overset{(\ref{AltBetakmo})}{=}& x\kk -\alpha\kk^{CG} g\kk + \alpha\kk^{CG} \frac{s\kk^TAg\kk}{s\kk^TAs\kk}s\kk,
\end{eqnarray}
which shows that at every iteration of CG, the search direction can be expressed as a linear combination of the negative gradient and the direction $s\kk$. It remains to show that the step sizes generated by the Flexible Scheme in \eqref{FlexSchemeCGsteps} are the same as used by CG in \eqref{xkCGgksk}. So,
\begin{eqnarray*}
    \alpha\kk^{CG} &\overset{(\ref{alphabeta})}{=}& -\frac{d\kk^T g\kk}{d\kk^TA d\kk}\\
    &\overset{(\ref{xdCG}+\ref{AltBetakmo})}{=}& -\frac{\left(-g\kk + \tfrac{s\kk^TAg\kk}{s\kk^TAs\kk}s\kk\right)^Tg\kk}{\left(-g\kk + \tfrac{s\kk^TAg\kk}{s\kk^TAs\kk}s\kk\right)^TA\left(-g\kk + \tfrac{s\kk^TAg\kk}{s\kk^TAs\kk}s\kk\right)}\\
    &=& - \frac{\tfrac{-g\kk^Tg\kk s\kk^TAs\kk + s\kk^TAg\kk s\kk^T g\kk}{(s\kk^TAs\kk)}}{\tfrac{g\kk^TAg\kk s\kk^T A s\kk - \left(s\kk^TAg\kk \right)^2  }{s\kk^TAs\kk}}\\
    &=& \frac1{\Delta}(g\kk^Tg\kk s\kk^TAs\kk - s\kk^T A g\kk s\kk^T g\kk)\\
    &\overset{\eqref{FlexSchemeCGsteps}}{=}& \alpha\kk^{(1)},
\end{eqnarray*}
as required.
To show the second step size, $-\alpha\kk \frac{s\kk^TAg\kk}{s\kk^TAs\kk} = \alpha\kk^{(2)}$ first note that by definition, the gradient is orthogonal to the previous search directions, so $d\kmo^Tg\kk =0$. Combining this with \eqref{AltDkmo}, and noting that $\alpha\kmo \neq 0$ gives
\begin{eqnarray}
 s\kk^Tg\kk = 0\label{CgOrthogonal}.
\end{eqnarray}
Now,
\begin{eqnarray*}
    -\alpha\kk^{CG} \frac{s\kk^TAg\kk}{s\kk^TAs\kk} &=& -\frac{g\kk^Tg\kk s\kk^TAs\kk - s\kk^T A g\kk s\kk^T g\kk}{\Delta} \frac{s\kk^TAg\kk}{s\kk^TAs\kk}\notag\\
    &\overset{(\ref{CgOrthogonal})}{=}& -\frac{g\kk^Tg\kk s\kk^TAs\kk}{\Delta}\frac{s\kk^TAg\kk}{s\kk^TAs\kk}\notag\\
    &=& \frac1{\Delta}(-g\kk^Tg\kk s\kk^TAg\kk)\notag\\
    &\overset{\eqref{CgOrthogonal}}{=}& \frac1{\Delta}(-g\kk^T g\kk s\kk^T A g\kk + g\kk^TAg\kk s\kk^T g\kk)\notag\\
    &\overset{\eqref{FlexSchemeCGsteps}}{=}& \alpha\kk^{(2)}.
\end{eqnarray*}
Hence, when $\ell = 0$ and $W\kk = \mat{g\kk,&s\kk}$ are used in the Flexible Scheme, CG is recovered. In a similar way, it can be shown that when $\ell = 1/2$ and $W\kk = \mat{g\kk,&s\kk}$, the Flexible Scheme recovers the conjugate residual method (see, e.g., \cite{Saad2003}).

While the theoretical results in this work guarantee convergence of CG, it is well known that one can obtain stronger theoretical results (i.e., a better rate of convergence) by taking into consideration the special properties of the search directions, see e.g., Theorem 7.8.3 in \cite{watkins2004fundamentals}.

Additionally, previous work has been done that extends the conjugate gradient methods to more general norms, including the $A^{2\ell-1}$-norms discussed in this work; see \cite{lim2024conjugate, faber1984necessary, ashby1990taxonomy, broyden1996new, fasano2015framework} for further details.

\subsubsection{Gradient Methods with Momentum and Acceleration}

Polyak's heavy-ball/momentum method \cite{polyak1964some}, and Nesterov's Accelerated Gradient Method (NAGM) \cite{Nesterov1983} are both popular algorithms that can be used to tackle problem \eqref{quadprob}. While neither of them fit the Flexible Scheme in Definition~\ref{def:Algorithms}, they are, nevertheless, worthy of mention here.

For the quadratic function \eqref{quadprob}, Nesterov's Accelerated Gradient Method can be written as 
\begin{align}
    z\kpo &= x\kk - \tfrac{1}{L}g(x\kk)\label{zkexpress}\\
    x\kpo &= z\kpo + \beta (z\kpo -z\kk)\label{xkexpress},
\end{align}
where $L = \lambda_1$ is the Lipschitz constant, and $\beta = \frac{\sqrt{\lambda_1} - \sqrt{\lambda_n}}{\sqrt{\lambda_1} +\sqrt{\lambda_n}}$ is an acceleration parameter. By substituting the expression for $z\kpo$ in \eqref{zkexpress}, into the expression for $x\kpo$ in \eqref{xkexpress}, the two iterative sequences can be combined into one as follows,
\begin{equation}\label{NAGM_1sequence}
    x\kpo = x\kk  - \tfrac1L g(x\kk) + \beta(x\kk - x\kmo) + \tfrac{\beta}L(g(x\kmo) - g(x\kk)).
\end{equation}
Through the lens of \eqref{NAGM_1sequence}, NAGM can be viewed as a multi-directional search method with the search directions $-g(x\kk)$, $x\kk - x\kmo$ and $g(x\kmo) - g(x\kk)$. The step lengths in \eqref{NAGM_1sequence} are constants, and as such, they do not satisfy condition \eqref{mina2lnorm} in Lemma~\ref{stepsizeminimizegrad} (i.e., the gradient is not minimized at the new point w.r.t. any matrix norm, as is required by Lemma~\ref{stepsizeminimizegrad} and Definition~\ref{def:Algorithms}), so NAGM does not fit the Flexible Scheme.
However, it can be shown that substituting the same search directions as used in NAGM, namely $W\kk = \mat{g\kk,& s\kk,& y\kk}$, where $s\kk = x\kk - x\kmo$ and $y\kk = g\kk - g\kmo$, into Definition~\ref{def:Algorithms} gives the conjugate gradient method. 

In a similar vein, substituting the same search directions as used in Polyak's heavy ball/momentum method into Definition~\ref{def:Algorithms} also recovers conjugate gradients.

\subsection{Forsythe $s$-gradient methods}
Forsythe's $s$--gradient methods \cite{Forsythe1968} use a search direction that is a linear combination of the sub-directions $g\kk,Ag\kk,\dots,A^{s-1}g\kk$, for some user specified parameter $s$. At each iteration, the step sizes are generated so as to minimize the function value over the hyper-plane spanned by the sub-directions, and as such, Forsythe's $s$--gradient methods belong to the class of algorithms given in Definition~\ref{def:Algorithms}.

As an example (with $\Pre = I$), choosing $s=3$ gives $W\kk = \mat{g\kk,& Ag\kk,& A^2g\kk}$, and using this in \eqref{GeneralisedStepsize} shows that the corresponding step sizes are given by 
\begin{equation}\label{Forsythesstepsizes}
    a\kk = \mat{g\kk^TA^{2\ell+1}g\kk & g\kk^TA^{2\ell + 2}g\kk & g\kk^TA^{2\ell + 3}g\kk\\
    g\kk^TA^{2\ell+2}g\kk & g\kk^TA^{2\ell + 3}g\kk & g\kk^TA^{2\ell + 4}g\kk\\
    g\kk^TA^{2\ell+3}g\kk & g\kk^TA^{2\ell + 4}g\kk & g\kk^TA^{2\ell + 5}g\kk}^{-1} \mat{g\kk^TA^{2\ell}g\kk \\ g\kk^TA^{2\ell+1}g\kk \\ g\kk^TA^{2\ell + 2}g\kk}.
\end{equation}
As is clear from \eqref{Forsythesstepsizes}, Forsythe's $s$--gradient methods require several matrix-vector products at every iteration, which is unfavourable, so they are not commonly used in practice. Further, these methods exhibit the same `zig-zagging' asymptotic behaviour as the $\ell-$MGD methods described in \cite{macdonald24} (see also Section~\ref{sec:relaxedellMG}), where the normalised gradients converge to the subspace spanned by the eigenvectors associated with the largest and smallest eigenvalues of $A$ \cite{Forsythe1968}, which leads to slow convergence \cite{Pronzato2009}.
It can be shown that using relaxation ($1\neq \omega \in (0,2)$) can lead to improvements in the practical behaviour of these methods (see, for example, Figure~\ref{fig:Forsythe2norm} in Section~\ref{sec:NumericalExp}). Finally, the works \cite{Pronzato2009,luenberger2008linear, chronopoulos1989s} establish that Forsythe's $s$-gradient method is equivalent to restarted CG.

\subsection{Extensions to General Norms}\label{sec:generalnorms}

The $A^{2\ell}$-norms are used in this work because they are concrete, computationally practical, and their use covers many of the well known algorithms in the literature. However, consider the following norm, defined in terms of a Laurent series, which is used to define the $P$-gradient algorithms in \cite{Pronzato_2005}.\footnote{There is a notation clash here, with $P$ being used to denote a preconditioner in this work, and $P$ defining a real function and associated family of algorithms in \cite{Pronzato_2005}, which we hope does not lead to confusion.}
\begin{definition}[Definition 1 in \cite{Pronzato_2005}]\label{def:Proznorm}
   Let $P(\cdot)$ be a real function defined on $[m, M]$, infinitely differentiable, with Laurent series
     $P(z) = \sum_{-\infty}^{\infty} c\kk z^k, c\kk \in \R \text{ for all } k,$
   such that $0 < \sum_{-\infty}^{\infty} c\kk a^k < \infty $ for $a \in [m, M]$. The $k$-th iteration of a $P$-gradient algorithm is defined by $x\kpo = x\kk -\gamma_k g\kk$
   where the step-length $\gamma\kk$ minimizes $\|g\kpo\|_{P(A)}^2$ with respect to $\gamma$, with $g\kpo =  \nabla f (x\kk - \gamma g\kk)$.
\end{definition}
Hilbert spaces are considered in \cite{Pronzato_2005}, but the results can be translated into $\R^n$ by taking $m = \lambda_n$, $M = \lambda_1$ and then choosing $P(A) = A^{2\ell-1}$; this recovers the $\ell$-MGD algorithms in \cite{macdonald24}, while $P(A) = \aln$ recovers the preconditioned variants in this work. Note that $P(A)$ in Definition~\ref{def:Proznorm} is positive definite, so it induces a norm, and it is shown in \cite[Appendix~A]{macdonald24} that the convergence theory for the relaxed $\ell$-MGD algorithms also holds under the more general $P(A)$ norms. 

Convergence of the Flexible Algorithmic Scheme in Definition~\ref{def:Algorithms} is established by showing that any algorithm fitting the scheme performs no worse than the corresponding preconditioned relaxed $\ell$-MGD method, and thus, the Flexible Scheme inherits the convergence guarantees of the preconditioned relaxed $\ell$-MGD method too. Consequently, convergence theory holds for any algorithm that fits the Flexible Scheme under the general $P(A)$ norm, Definition~\ref{def:Proznorm} (not just the $A^{2\ell-1}$-norm).

\subsection{Generalised Delayed Weighted Gradient Method}
The Generalised Delayed Weighted Gradient Method (GDWGM) \cite{oviedo2022family, andreani2022extended} is a method that exhibits finite termination. The method is constructed by taking an exact stepsize in the gradient direction, then using an over-relaxation scheme. Each of the coefficients is found by minimizing a merit function $F_\mu = (1-\mu)E(x) + \mu \|\nabla f(x)\|_2^2$, where $E(x) = \tfrac12 (x-x^*)^TA(x-x^*)$ and $\mu \in [0,1]$. The steps of the algorithm are
\begin{align*}
    z\kk &= x\kk - \alpha\kk \nabla f(x\kk)\\
    x\kpo &= \beta\kk z\kk + (1-\beta\kk)x\kmo,\\
\end{align*}
where $\alpha\kk = \arg \min\limits_{\alpha > 0} F_\mu(x\kk - \alpha \nabla f(x\kk))$ and $\beta\kk = \arg \min\limits_{\beta} F_\mu (\beta z\kk + (1-\beta)x\kmo)$.

Note that varying $\mu \in [0,1]$ gives a family of methods, and the specific choice $\mu=0$ recovers conjugate gradients.
The merit function can be, equivalently, described in terms of norms as 
\begin{eqnarray*}
    2F_\mu &=& (1 - \mu)2E(x) + 2\mu \|\nabla f(x\kk)\|^2_2 \\
    &=& (1-\mu) \|x\kk -x^*\|_A^2 + 2\mu \| \nabla f(x\kk) \|^2_2 \\
    &=& (1-\mu) \| g\kk \|_{A^{-1}}^2 + 2\mu \| g\kk \|^2_2\\
    &=& \|g\kk\|^2_{A^{-1}W_\mu}, 
\end{eqnarray*}
where $W_\mu = (1-\mu)I + 2\mu A$, as defined in \cite{oviedo2022family}. This norm is a special case of the more general $P(A)$ norms (let $P(A) = (1-\mu)A^{-1} + 2\mu I = A^{-1}W_\mu$) induced through a Laurent series (see Section~\ref{sec:generalnorms}). Substituting $W\kk = \mat{g\kk, & s\kk}$, where $s\kk = x\kk - x\kmo$ into the Flexible Scheme in Definition~\ref{def:Algorithms}, and minimizing under the $A^{-1}W_\mu$ norm of the gradient gives the step sizes (without preconditioning $P = I$) 
\begin{equation}
    a\kk = \mat{g\kk^TW_\mu A g\kk & s\kk^T W_\mu A g\kk\\
    s\kk^TW_\mu A g\kk & s\kk^T W_\mu A s\kk}^{-1} \mat{g\kk^T W_\mu g\kk \\ s\kk^T W_\mu g\kk},
\end{equation} 
and thus, the GDWGM \cite{oviedo2022family, andreani2022extended} is recovered. Direct computation of the step sizes, as done in Section~\ref{sec:CG}, shows equivalence between our scheme and GDWGM.

The paper \cite{ashby1990taxonomy} generalises the conjugate gradient method to a wider family of norms. The Pronzato norms are contained in this wider family, and the $A^{-1}W_\mu$ norms that give the GDWGM methods are a specific case of the Pronzato $P(A)$ norms. The Flexible Algorithmic Framework shows that if the search directions are chosen to be $W\kk = \mat{g\kk, & s\kk}$, it results in a conjugate gradient method (with a different norm). This means we get stronger convergence results (see e.g., Theorem 7.8.3 in \cite{watkins2004fundamentals}) and finite termination in exact arithmetic (Theorem 7.8.12 in \cite{watkins2004fundamentals}). For the GDWGM method, this means our scheme directly shows it has a stronger convergence rate and finite termination.

\subsection{New algorithms fitting the Flexible Algorithmic Framework}
\label{section:newmethods}

It is widely acknowledged that Conjugate Gradients is the algorithm of choice for strongly convex quadratic minimization. Nevertheless, when considering novel algorithms for more general functions, during the development phase, investigating the performance of said algorithms on quadratics can provide valuable insights into practical algorithmic behaviour. 

This work provides a useful tool for algorithm development. Any algorithm that a user can imagine, which fits the Flexible Scheme in Definition~\ref{def:Algorithms}, is automatically supported by convergence guarantees and a linear rate of convergence on quadratic functions, so the user can simply push ahead with the task of numerical testing. Thus, this work can assist with streamlining algorithm development.

To demonstration this, three prototype algorithms are considered now. All algorithms fit the Flexible Scheme in Definition~\ref{def:Algorithms}, so they are equipped with the theoretical convergence guarantees already established in this work. They are described here, and preliminary numerical experiments for them are presented in Section~\ref{sec:NumericalNew}.

\paragraph{Gradient Descent with a random direction.} Gradient Descent is known to be slow near saddle points, and Perturbed Gradient Descent \cite{jin2017escape} was developed to try to help GD escape saddle points more rapidly. At every iteration, the algorithm takes the negative gradient, and adjusts it by adding to it a random direction, and then a step is taken along this (single) perturbed gradient direction. Therefore, Perturbed Gradient Descent does not fit the Flexible Scheme in Definition~\ref{def:Algorithms} because $g\kk \not \in {\rm col} (W\kk)$ (i.e., Assumption~\ref{Assumepi} fails). Nevertheless, it does inspire the following, distinct but related algorithm, namely, Gradient Descent with a Random Direction (GD--RD). 

At every iteration of GD--RD, a random direction $r\kk$ is sampled, and $W\kk$ is taken to be $W\kk = \mat{g\kk, &r\kk}$, so Assumption~\ref{Assumepi} is satisfied. For concreteness, setting $P = I$, and taking any fixed $\ell \in \{0,\tfrac12,1,\tfrac32,\dots\}$, the step sizes are computed from \eqref{GeneralisedStepsize} to be
\begin{equation}
    a\kk = \mat{g\kk^T A^{2\ell+1}g\kk & g\kk^TA^{2\ell+1} r\kk \\
    r\kk^TA^{2\ell + 1}g\kk & r\kk^TA^{2\ell + 1}r\kk}^{-1} \mat{g\kk^T A^{2\ell} g\kk \\ r\kk^T A^{2\ell} g\kk}.
\end{equation}
That is, at every iteration of GD--RD, two sub-search directions are used ($g\kk$ and $r\kk$), and the step sizes are generated so that the $A^{2\ell-1}$-norm of the gradient is minimized over the two-dimensional subspace spanned by $g\kk$ and $r\kk$. (This contrasts with Perturbed Gradient Descent, which takes a step along the one dimensional subspace spanned by the single perturbed gradient direction $g\kk + r\kk$.) Preliminary numerical experiments are presented for GD--RD in Section~\ref{sec:NumericalNew}.

\paragraph{Forsythe's $s=2$ method with momentum.} A second algorithm considered now is motivated by the positive aspects of Forsythe's method, coupled with the benefits of momentum. For this new algorithm, hereby referred to as Forsythe's $s=2$ method with momentum, at every iteration the sub-search directions correspond to $W\kk = \mat{g\kk, &Ag\kk, &x\kk-x\kmo}$. Figure~\ref{fig:tripdirection} in Section~\ref{sec:NumericalNew} demonstrates the practical behaviour of this algorithm.

\paragraph{A momentum type method with a random direction.} Lastly, a third prototype algorithm considered here involves taking a momentum type algorithm, and adding in a random direction. Specifically, the momentum type-random direction method considered here sets $W\kk = \mat{g\kk, &x\kk-x\kmo, &r\kk}$ at every iteration. Numerical experiments for this algorithm are presented in Section~\ref{sec:NumericalNew}.

\section{Numerical Experiments}\label{sec:NumericalExp}
Extensive numerical experiments can be found in the literature for the existing state-of-the-art methods mentioned previously (see, for example, any of the references listed in this work). The purpose of the numerical experiments presented now is to give a (deliberately brief) illustration of the typical numerical behaviour of some of the methods, and provide intuition as to how multiple sub-search directions can impact algorithm performance. All the code is written in Python 3.7.4 using an AMD Ryzen 5 3600 CPU with 16GB of RAM.

The numerical set up is as follows. Let $m=1200$ and $n=1000$, let $B \in \R^{m \times n}$ be a random matrix with uniformly distributed entries in $[0, 1)$, and then form the symmetric positive definite matrix $A = B^TB$. The minimizer $x^*\in \R^n$ is generated from the same distribution, and then $b = Ax^*$ is computed. All algorithms are given the same randomly generated initial vector $x\kz\in \R^n$, and for each problem instance/data set a new random point is initialized. The algorithms run until a tolerance of $\|g\kk\|_2^2 < \epsilon = 10^{-6}$ is reached. If an algorithm did not reach the stopping tolerance within 1000 iterations, it was terminated. 
 
\subsection{Numerical experiments comparing existing algorithms}

\begin{figure}[h!]
    \centering
    \includegraphics[width=0.48\linewidth]{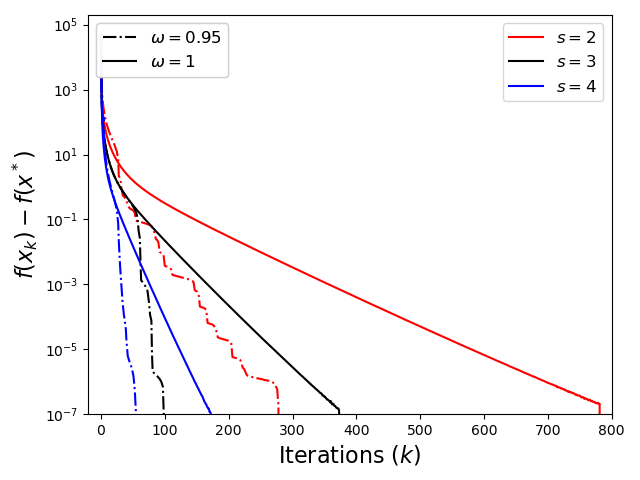}
    \includegraphics[width=0.48\linewidth]{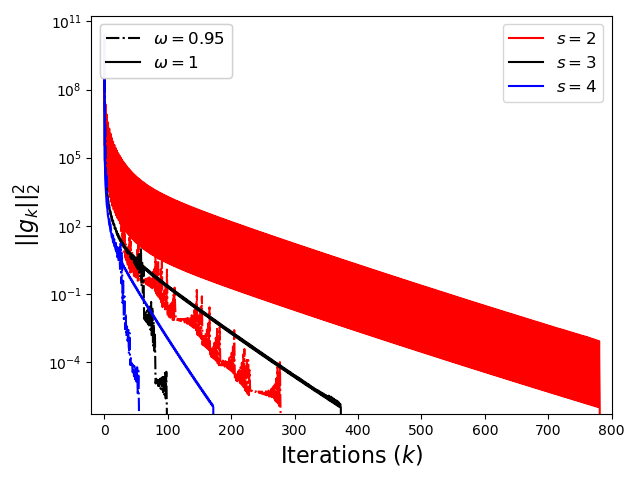}
    \caption{Forsythes $s-$gradient method with relaxation $\omega = 0.95$, compared with the unrelaxed methods. Function values are shown on the left, and the (squared) gradient norm is shown on the right.}
    \label{fig:Forsythe2norm}
\end{figure}

\paragraph{Forsythe's $s$--gradient methods with relaxation.}
This work extends the results of \cite{Forsythe1968}, and shows that Forsythe's $s$--gradient methods are guaranteed to converge when fixed relaxation $\omega \in (0,2)$ is employed. Further, the numerical experiment presented here shows that Forsythe's $s$--gradient methods can indeed benefit from relaxation in practice. Consider the following experimental set up. Set $\ell = 0$, and consider Forsythe's $s=2$--, $s=3$-- and $s=4$--gradient methods (defined by setting $W\kk = \mat{g\kk,& Ag\kk,& \dots,& A^{s-1}g\kk}$ at every iteration, for $s=2$, $s=3$ and $s=4$, respectively). Each of the three algorithms is compared with no relaxation ($\omega = 1$) and with fixed relaxation of $\omega = 0.95$ (giving a total of 6 different algorithm instances fitting the Flexible Scheme in Definition~\ref{def:Algorithms}). The results are shown in Figure~\ref{fig:Forsythe2norm}.

Figure~\ref{fig:Forsythe2norm} shows the evolution of $f(x\kk)-f(x^*)$ as iterations progress (left plot) and the evolution of $\|g\kk\|_2^2$ as iterations progress (right plot). The advantages of using an enriched search direction is clearly seen, with the $s=4$ algorithm variant (solid blue line) (using the 4 sub-search directions $g\kk$, $Ag\kk$, $A^2g\kk$ and $A^3g\kk$) reaching the stopping threshold in many fewer iterations than both the $s=3$ (solid black line) and $s=4$ (solid red line) variants. (Of course, there is a trade-off, with the $s=4$ variant requiring more matrix vector products than both the $s=3$ and $s=2$ algorithms.)

Moreover, Figure~\ref{fig:Forsythe2norm} also highlights the benefits of relaxation, with each of the relaxed ($\omega=0.95$) algorithm variants (the dash-dot lines) requiring fewer iterations to converge compared with the un-relaxed ($\omega=1)$ variants (solid lines).

Lastly, it is well known that the norm of the gradient can increase for Forsythe's $s=2$ method, and note that this erratic behaviour is captured in the right plot in Figure~\ref{fig:Forsythe2norm}, where the `red-band' is actually a single solid line that decreases non-monotonically as iterations progress.

\paragraph{A comparison of conjugate direction methods.}

The following experiment compares the performance of several conjugate direction methods. In particular, setting $W\kk=\mat{g\kk,&s\kk}$ in Definition~\ref{def:Algorithms} with either $\ell$ fixed at 0 or 1/2 gives Conjugate Gradients and Conjugate Residuals, respectively. The $\ell=1$ Conjugate Direction (CD) method was also considered, as was GDWGM from \cite{oviedo2022family} with the choice $\mu = 0.5$. Each of the four algorithms was compared either without a preconditioner ($P=I$), or with the Jacobi preconditioner (i.e., $P$ is a diagonal matrix with diagonal entries $P_{ii} = A_{ii}$).

\begin{figure}[h!]
    \centering
    \includegraphics[width=0.48\linewidth]{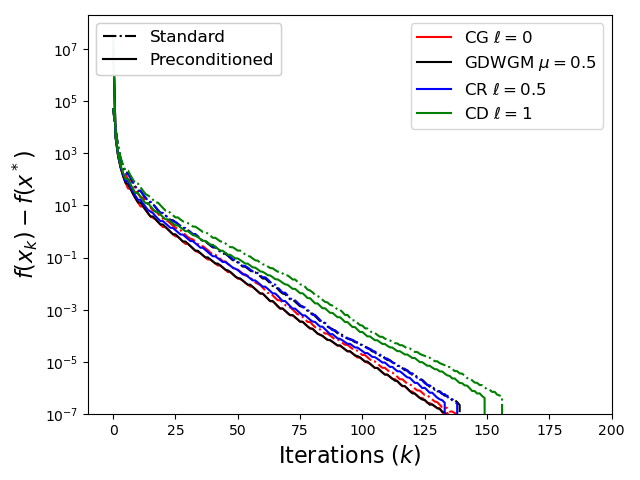}
    \includegraphics[width=0.48\linewidth]{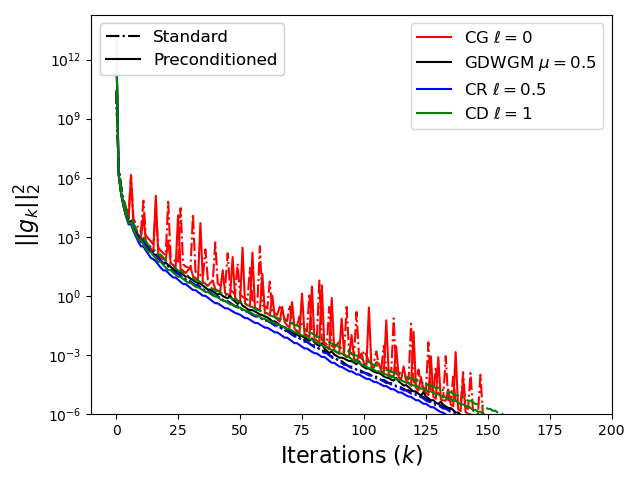}
    \caption{Conjugate direction methods, both with and without preconditioning. The function values are shown to take slightly fewer iterations when using a Jacobi preconditioner.}
    \label{fig:CDPrecon}
\end{figure}

The results are shown in Figure~\ref{fig:CDPrecon} where the left plot shows the evolution of $f(x\kk)-f(x^*)$ as iterations progress and the right plot shows the evolution of $\|g\kk\|_2^2$ as iterations progress. All four of these conjugate direction methods takes a similar number of iterations to converge, with only a modest further reduction in iterations for the preconditioned variants (solid lines).

\subsection{Numerical experiments comparing new algorithms}
\label{sec:NumericalNew}

\paragraph{Comparing the relaxed $\ell$-MGD methods and GD with a random direction.}

\begin{figure}[h!]
    \centering
    \includegraphics[width=0.48\textwidth]{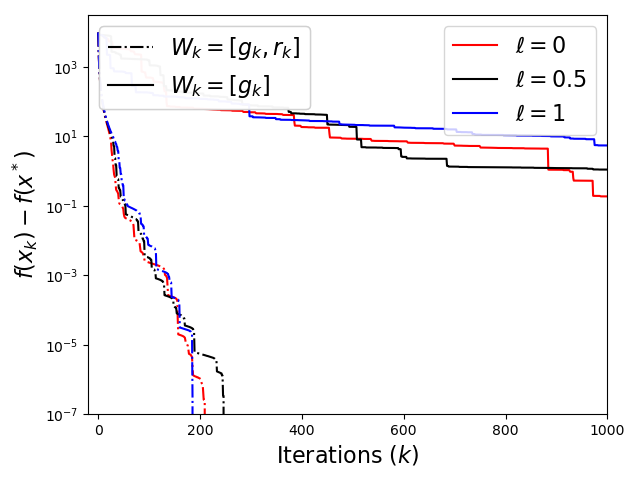}
    \includegraphics[width=0.48\textwidth]{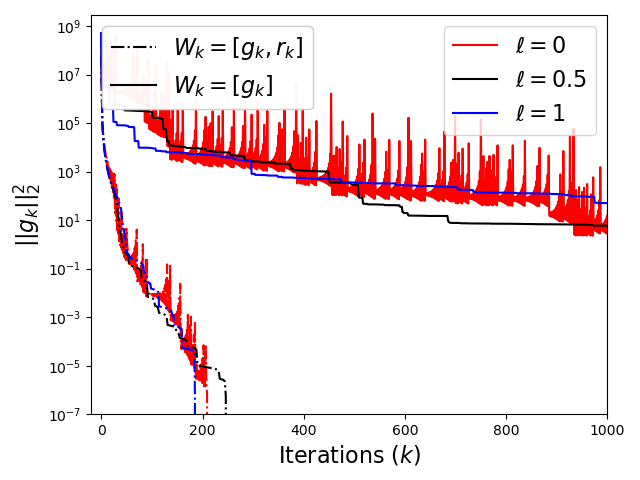}
    \caption{A comparison of several relaxed $\ell$-MGD methods, with the Gradient Descent with a random direction method. The left plot shows the evolution of the function values, and the right plot shows the evolution of the squared 2-norm of the gradient as iterates progress. }
    \label{fig:randExample}
\end{figure}

The goal of this experiment is to investigate whether the performance of gradient based methods can be improved when randomization is used. In particular, here, the behaviour of several relaxed $\ell$-MGD algorithms (see Section~\ref{sec:relaxedellMG} and \cite{macdonald24}), which take $W\kk = g\kk$ at every iteration, are compared with the GD with a random direction method (recall Section~\ref{section:newmethods}) which takes $W\kk = \mat{g\kk,&r\kk}$ at every iteration. For this experiment, the entries of $r\kk\in\R^n$ are sampled from a normal distribution at every iteration. Each choice of $W\kk$ was used in the Flexible Scheme in Definition~\ref{def:Algorithms} with $\ell=0$, 1/2 and 1, giving a total of 6 algorithm variants to compare. In all cases, relaxation of $\omega = 0.95$ was employed, and the results are shown in Figure~\ref{fig:randExample}.

Figure~\ref{fig:randExample} shows the evolution of the function values (left plot) and evolution of the squared 2-norm of the gradient (right plot) as iterations progress, for the algorithms just described. It is clear from Figure~\ref{fig:randExample} that including a second random direction in a gradient based method (i.e., using GD-RD) is beneficial, with all of the GD-RD variants (for $\ell=0$, 1/2 and 1) meeting the stopping condition in approximately 200 iterations, whereas the relaxed $\ell$-MGD algorithms fail to make meaningful progress in the maximum allowed iterations. It should be noted that, while GD-RD requires an additional matrix-vector product at every iteration compared with the corresponding relaxed $\ell$-MGD algorithm, this is compensated for by the significant drop in the number of iterations required to meet the convergence threshold.

\begin{figure}[h!]
    \centering
    \includegraphics[width=0.48\linewidth]{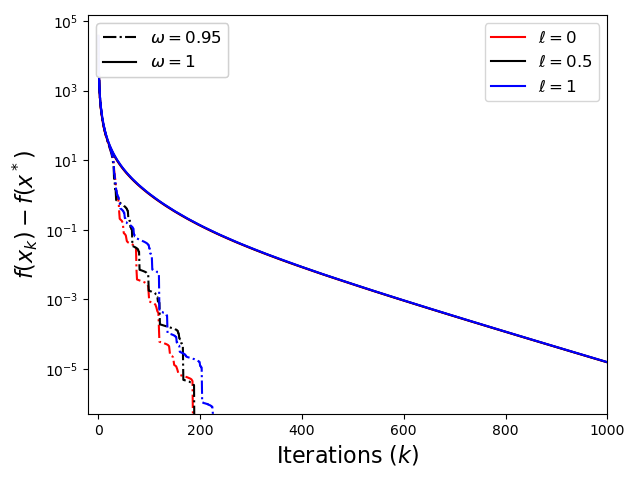}
    \includegraphics[width=0.48\linewidth]{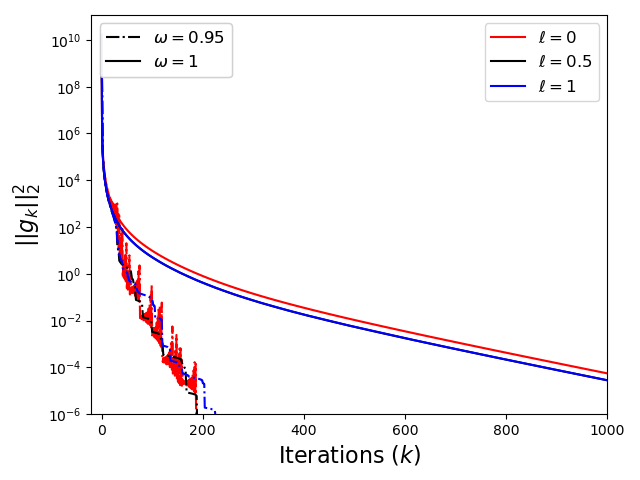}
    \caption{Results of perturbed gradient descent with $(\omega = 0.95)$ and without relaxation $(\omega = 1)$. Function values shown on the left and gradient norm is shown on the right.}
    \label{fig:RelaxedRandomDirec}
\end{figure}

For completeness, Figure~\ref{fig:RelaxedRandomDirec} shows the performance of GD-RD algorithm, for each of $\ell = 0, 1, 1/2$, both with $(\omega = 0.95)$ and without $(\omega = 1)$ relaxation. Figure~\ref{fig:RelaxedRandomDirec} shows the benefits of relaxation, with the un-relaxed variants requiring more than the maximum allowed (1000) iterations to reach the stopping tolerance.

\paragraph{Comparison of Forsythe's $s=2$ method with momentum and a momentum type method with a random direction.}

\begin{figure}[h!]
    \centering
    \includegraphics[width=0.48\linewidth]{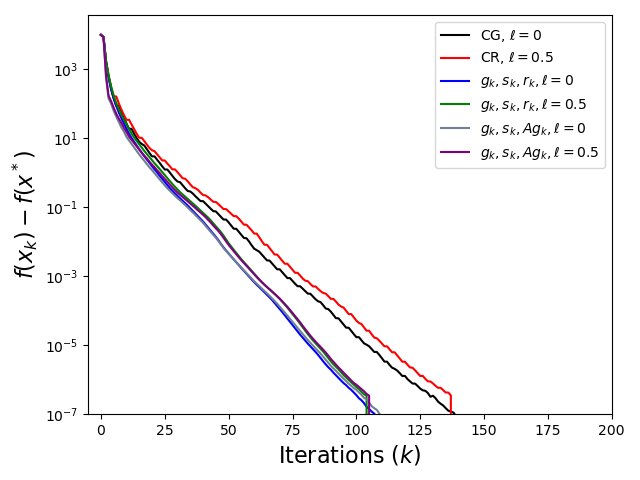}
    \includegraphics[width=0.48\linewidth]{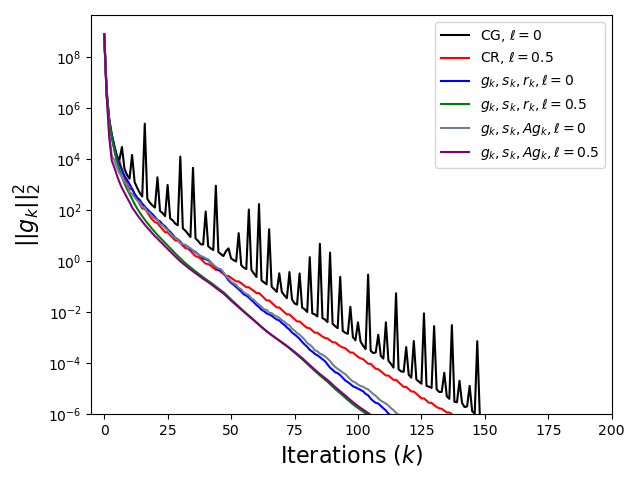}
    \caption{The CG and CR methods, compared with the `enhanced' methods, which take $Ag\kk$ or $r\kk$ as an extra direction. The decrease in function values is shown on the left, and the gradient norm is shown on the right.}
    \label{fig:tripdirection}
\end{figure}

Forsythe's $s=2$ method with momentum takes $W\kk = \mat{g\kk, &Ag\kk, &x\kk-x\kmo}$ at every iteration, while a momentum type method with a random direction sets $W\kk = \mat{g\kk, &x\kk-x\kmo, &r\kk}$. Each choice of $W\kk$ is combined with $\ell=0$ and $\ell = 0.5$, giving 4 different algorithms fitting the Flexible Scheme. As a baseline (and recalling the discussion in Sections~\ref{sec:CG}), these methods are compared with CG (which fixes $\ell = 0$ and uses $W\kk = \mat{g\kk &x\kk-x\kmo}$ at every iteration) and CR (which fixes $\ell = 1/2$ and uses $W\kk = \mat{g\kk &x\kk-x\kmo}$ at every iteration). These 6 algorithms are compared, and the results are shown in Figure~\ref{fig:tripdirection}.

The left plot in Figure~\ref{fig:tripdirection} shows the evolution of $f(x\kk)-f(x^*)$ as iterations progress. It also shows that adding either $r\kk$ or $Ag\kk$ to the `CG directions', $g\kk$ and $s\kk$, (under either $\ell=0$ or $\ell=1/2$), is leads to a reduction in the number of iterations needed for convergence, compared with CG or CR. (Figure~\ref{fig:tripdirection} shows that conjugate gradient methods take around 130 iterations to converge, whereas the new methods take around 100 iterations to converge.) However, it should be noted that the new methods require more matrix-vector products per iteration, so they have a higher computational cost per iteration compared with CG and CR.

\section{Conclusion}
\label{sec:Conclusion}

This work considers the minimization of strictly convex quadratic functions and presents a flexible algorithmic framework that can be used to solve such problems. The framework allows the user to choose the number of sub-search directions, as well as which directions to use (both of which can change from one iteration to the next), and the scheme also allows relaxation and preconditioning to be employed. Thus, the scheme is general and flexible. Moreover, as long as the gradient is included as one of the sub-search direction at every iteration, any algorithm that fits the framework is equipped with theoretical guarantees that show convergence to the minimizer at a linear rate. It is also shown that several existing state-of-the-art methods fits the framework (for example, steepest descent and the conjugate gradients method), and several new algorithms are also proposed that fit into the scheme. Our results provide for an efficient and robust approach to the analysis of gradient based algorithms: when a new algorithm is devised for a specific or general class of applications, checking its adherence to the Flexible Scheme is sufficient for guaranteed convergence, and efforts can be directed to testing and refining detailed parameter choices for best performance.

\bibliographystyle{unsrt}
\bibliography{refs}

\end{document}